# FOUR POINTS AND A QUADRIC

Nicholas Phat Nguyen

**Abstract.** This paper provides a simple proof of Hesse's theorem in projective geometry for any dimension.

**1.   INTRODUCTION.** Almost anyone studying or teaching projective plane geometry will have come across the following theorem, named after Ludwig Otto Hesse (1811-1874), a German mathematician whose name is associated with number of important mathematical concepts, such as the Hessian matrix of second-order derivatives of a scalar-valued function of several variables.

**Theorem (Hesse's Theorem for Quadrangle and Conic in Projective Plane):** *If two pairs of opposite sides of a quadrangle in a projective plane are conjugate lines relative to a conic, then the third pair of opposite sides is a also a pair of conjugate lines relative to the same conic.*

In this note, we will revisit this classic theorem and provide a simple proof of Hesse's theorem for a quadric in projective space of any dimension.   We hope that even readers with little knowledge of projective geometry could still follow the main ideas in this note, and enjoy making some connection between the old and the new.

Let us first consider what Hesse's theorem means.  Let *A, B, C* and *D* be four distinct points, i.e., a quadrangle, in a projective plane.  These four points give us three pairs of lines corresponding to the three partitions of the quadrangle into sets of two, namely the pairs (line *AB*, line *CD*), (line *AC*, line *BD*), and (line *AD*, line *BC*).  The two lines in each pair may not be distinct, depending on the configuration of the quadrangle.  Two lines $l$ and $m$ are said to be conjugate with respect to a conic *Q* if there is a point on the line $l$ that is orthogonal (relative to *Q*) to all the points of the line $m$, or equivalently if there is a point on the line $m$ that is orthogonal to all the points of the line $l$.  Hesse's theorem says that out of the three pairs of lines that we can form out of a quadrangle, if two pairs are pairs of conjugate lines with respect to a conic *Q*, then the remaining pair is also a pair of conjugate lines with respect to *Q*.

For example, suppose we take the quadrangle *A*, *B*, *C*, *H* in the real affine plane $\mathbf{R}^2$ regarded as a subset of the real projective plane $\mathbf{P}(\mathbf{R}^2)$, and suppose that *H* is the intersection of the altitude line from *B* to *AC* and the altitude line from *C* to *AB*. Now take the conic *Q* to be any circle centered at *H*. The line *BH* is a diameter of the circle *Q*, and therefore is conjugate to any line perpendicular to it, such as the line *AC*. Similarly, the diameter *CH* is conjugate to the line *AB*. Hesse's theorem then says that the remaining pair (line *AH*, line *BC*) must also be a pair of conjugate lines with respect to the circle *Q*, which means the diameter *AH* is perpendicular to the line *BC*. In other words, the three altitude lines of the triangle *ABC* must be concurrent at *H*. Accordingly, Hesse's theorem could be regarded as a generalization of the triangle altitude theorem in Euclidean plane geometry.

Hesse's theorem is often also stated in a dual form by replacing point with line and line with point. Specifically, in the dual version, we have four distinct lines in a projective plane (a quadrilateral), meeting in six points called vertices (which may not be all distinct). These six vertices can be grouped in three pairs of opposite vertices based on partitions of the lines in the quadrilateral. If two pairs of such opposite vertices are pair of orthogonal points relative to a given conic, then the remaining pair of points are also likewise orthogonal to each other.

In practice, most versions of Hesse's theorem in the literature either explicitly or implicitly require that the four given points be in general position, with no three points colinear (this configuration is known as a complete quadrangle), and that the conic is non-degenerate. A direct proof Hesse's theorem in the literature would typically require the drawing of an extra line, a line that is polar to one of the four given points. The argument would then apply the Desargues' involution theorem to the intersections of that extra line with the three pairs of opposite sides in our quadrangle configuration to deduce the required result. That approach certainly requires the conic to be non-degenerate (so that the polar of a point is a line), and the four given points to be in general position. These conditions and the hypothesis of the theorem would effectively preclude the possibility of all four given points being on the conic, where the Desargues involution argument would completely fail to apply. However, this is a fine point that takes some work to go through, so many proofs simply just ignore it.

Hesse's theorem applies not just to four points and a conic in a projective plane, but also to four points and a quadric surface in projective space of dimension 3. In the classic treatise *Principles of Geometry*, Henry Baker stated and proved such a general version of Hesse's theorem for four distinct points and a regular quadric surface in projective space of dimension 3. However, the proof in this higher dimensional case as outlined by Baker is considerably more difficult and elaborate, requiring by his own count a configuration of 15 points, 20 lines and 15 planes. See [1] at pp. 34-36.

There is actually a very simple underlying idea for Hesse's theorem in either dimension. Specifically, we will show that:

**Theorem (General Version of Hesse's Theorem):** *Consider a projective space of any dimension over a field of characteristic ≠ 2. Suppose we have four distinct points in that projective space and a quadric. If two pairs of opposite lines determined by these four points are conjugate lines relative to the quadric, then the third pair of opposite lines is a also a pair of conjugate lines.*

This general version of Hesse's theorem does not impose any condition on the quadric or on the four given points other than that they are distinct. If the reader feels some unease about the concept of conjugate lines, please do not worry. That concept will be clarified as part of our proof of Hesse's theorem below.

**2. PROOF OF HESSE'S THEOREM.** We will work over a ground field $K$ of characteristic ≠ 2. A projective space of dimension $n$ defined over $K$ is the set $P(E)$ of lines through the origin in a $K$-vector space $E$ of linear dimension $(n + 1)$.

A quadric in $P(E)$ is a subset of $P(E)$ consisting of all zeros of a homogeneous polynomial $q(x)$ of degree 2, where the variables represent the coordinates of each point of the vector space $E$ relative to some basis. Because the ground field $K$ has characteristic $\neq 2$, we can express $q(x)$ as $<x, x>$, where $(u, v) \mapsto <u, v>$ represents a symmetric bilinear form on $E$, i.e., a mapping from the product space $ExE$ into the ground field $K$ that is linear in each argument and symmetric.

Now suppose we have four distinct points $a, b, c, d$ in the projective space $P(E)$, together with a quadric $Q$ induced by a symmetric bilinear form $(u, v) \mapsto <u, v>$. We begin with a lemma that characterizes conjugate lines relative to the quadric $Q$.

*Lemma*: *The line ab contains a point that is orthogonal to all points on the line cd if and only if $<a, c><b, d> - <a, d><b, c> = 0$, where the expression is calculated by taking any representative vectors of the four points.*

The pairing of two points in the projective space is generally not well-defined, because different representative vectors for the same points will give us different pairing values. However, note that the expression $<a, c><b, d> - <a, d><b, c>$ has certain symmetry, so that a change in representative vectors will change each term in the expression by the same multiplicative factor. Therefore, if some representative vectors for these four points give us zero value when we apply the expression to them, then any other set of representative vectors will give us the same zero value for the expression.

*Proof of Lemma*: Let $A, B, C, D$ be a set of representative vectors for the points $a, b, c, d$. Any point in the line $ab$ can be represented by a vector of the form $\alpha A + \beta B$, where $\alpha$ and $\beta$ are elements of the ground field $K$ that are not both zero. For such a point to be orthogonal to all points on the line $cd$, it is necessary and sufficient that it is orthogonal to $c$ and to $d$, i.e., $<\alpha A + \beta B, C> = 0$ and $<\alpha A + \beta B, D> = 0$.

Expanding out these expressions, we have:

$<A, C>\alpha + <B, C>\beta = 0$

$<A, D>\alpha + <B, D>\beta = 0$

From linear algebra, we know that such a system of two linear equations in two variables have a non-trivial solution (where $\alpha$ and $\beta$ are not both zero) if and only if the determinant $<A, C><B, D> - <A, D><B, C> = 0$. ∎

By symmetry of the expression, if the line *ab* contains a point that is orthogonal to all points on the line *cd*, then by the lemma the line *cd* also contains a point that is orthogonal to all points on the line *ab*.  This is by itself a non-trivial fact.  We will refer to the lines *ab* and *cd* as conjugate lines relative to the quadric if this is the case.

We now return to our general configuration of four points *a, b, c, d* and a quadric *Q*, defined by symmetric bilinear form $(u, v) \mapsto <u, v>$.  If the line *ab* and the line *cd* are conjugate lines relative to *Q*, then by the above lemma we have:

$<a, c><b, d> - <a, d><b, c> = 0$, or $<a, c><b, d> = <a, d><b, c>$.

Similarly, if the line *ac* and the line *bd* are conjugate lines, then

$<a, b><c, d> - <a, d><b, c> = 0$, or $<a, b><c, d> = <a, d><b, c>$.

These two equations together imply that we must have

$<a, c><b, d> = <a, b><c, d>$.

That means the line *ad* and the line *bc* must also be conjugate lines relative to *Q*, and our proof is complete. ∎

3.  **HESSE'S THEOREM IN PROJECTIVE DIMENSION ONE**.  Note that despite the lack of restriction on dimension in the general version of Hesse's theorem, we really do not gain anything beyond projective dimension 3.  That is because four distinct points in a projective space will generate a projective subspace of at most dimension 3, and we can just confine our attention to that subspace, where all relevant things for the theorem take place.

Let's consider what happens in projective dimension one.  In that case, we are in a projective line, and all the lines in Hesse's theorem are one and the same line.  Therefore,

Hesse's theorem is automatically true in that situation. The ambient projective line is conjugate to itself relative to a quadric if and only if there is a point in the projective line orthogonal to everything, which is the case if and only if the pairing that defines the quadric is degenerate.

If we look at the projective line in question as the set of lines in a vector space of dimension 2, then the proof of Hesse's theorem gives us the following result for a quadratic space of dimension 2, which can be regarded as Hesse's theorem for a projective line.

*Proposition*: *Let E be a vector space of dimension two over a ground field K of characteristic ≠ 2, endowed with a symmetric bilinear form Q(u, v). The form Q is degenerate if and only if for any four vectors a, b, c, d in E, no two of which are proportional, we have Q(a, c)Q(b, d) = Q(a, d)Q(b, c).*

Note that the usual criterion for $Q$ to be degenerate is that the determinant of its matrix relative to any basis of $E$ be zero. The above proposition, a byproduct of the general version of Hesse's theorem, gives us another criterion. This is an elegant and useful result, but it does not seem to be noted in the literature on quadratic forms.

Now suppose $E = K^2$ and consider the standard scalar product $S(x, y) = x_1.y_1 + x_2.y_2$ on $E$. This scalar product is a non-degenerate symmetric bilinear form. If we take $u = (1, 0)$ and $v = (0, 1)$ as the standard basis for E, then for vectors $x, y$ such that the four vectors $u, v, x, y$ are projectively distinct (no two of which are proportional to each other), we have

$S(u, x).S(v, y) = x_1.y_2$ , and

$S(u, y).S(v, x) = x_2.y_1$

Because $S$ is non-degenerate, the proposition above implies that $x_1.y_2 \neq x_2.y_1$. Note that the cross-ratio of the four projective points represented by $u, v, x, y$ in that order is simply the ratio $(x_2.y_1)/(x_1.y_2)$. So what Hesse's theorem in dimension one means in this situation is simply that the cross ratio of any four distinct points in a projective line is never equal to one. But that is of course a basic property of the cross ratio (a cross ratio can take on any value except for 1, 0 and ∞).